\documentclass[dvips]{arxbj}
\usepackage{graphics}

\aid{0}
\volume{16}
\issue{2}
\pubyear{2010}
\firstpage{435}
\lastpage{458}
\doi{10.3150/09-BEJ217}

\makeatletter

\theoremstyle{definition}

\makeatother

\begin{document}
\begin{frontmatter}

\title{Latent diffusion models for survival analysis}
\runtitle{Latent diffusion models for survival analysis}

\begin{aug}
\author[1]{\fnms{Gareth O.} \snm{Roberts}\thanksref{1}\ead[label=e1]{Gareth.O.Roberts@warwick.ac.uk}}
\and
\author[2]{\fnms{Laura M.} \snm{Sangalli}\thanksref{2}\ead[label=e2]{laura.sangalli@polimi.it}\corref{}}
\runauthor{G.O. Roberts and L.M. Sangalli}
\address[1]{CRiSM, Department of
Statistics, University of Warwick, Coventry CV4 7AL,
UK.\\ \printead{e1}}
\address[2]{MOX, Dipartimento di
Matematica, Politecnico di Milano, P.zza L. da Vinci 32, 20133
Milano, Italy.  \printead{e2}}
\end{aug}

\received{\smonth{6} \syear{2008}}
\revised{\smonth{2} \syear{2009}}

%
\begin{abstract}
We consider Bayesian hierarchical models for
survival analysis, where the survival times are modeled through an
underlying diffusion process which determines the hazard rate. We
show how these models can be efficiently treated by means of Markov
chain Monte Carlo
techniques.
\end{abstract}

%

\end{frontmatter}

\section{Introduction}
\label{sec:intro}

Diffusion processes have found many applications in the modeling of
continuous-time phenomena, for problems related to a variety of
scientific areas, ranging from economics to biology, from physics to
engineering. Here, we use diffusion processes as building blocks for
the definition of models for survival and event history analysis.
This idea is not new (see, e.g., the reviews in
\citeauthor{Aalen-Gjessing-2001} (\citeyear{Aalen-Gjessing-2001},
\citeyear{Aalen-Gjessing-2004})). However, in this
paper, we are able to considerably extend the flexibility of the
diffusion models used, by adopting powerful Markov chain Monte Carlo
techniques.

Diffusion models for survival analysis have been proposed because,
as summarized in \citet{Aalen-Gjessing-2004}, ``when modelling
survival data it may be of interest to imagine an underlying process
leading up to the event in question.''
Such a process might, for example,
represent the development of a disease. Two types of models have
been considered in the literature: models where the event happens
when a diffusion process hits some barrier and models where the
hazard rate is some suitable function of the diffusion. For the
former type of model, we refer the reader to
\citet{Aalen-Gjessing-2001}, \citet
{Aalen-Borgan-Gjessing-2008} and
references therein. Here, we are interested in the latter.
\citet{Woodbury-Manton-1977} proposed a model where the hazard rate
is a quadratic function of an Ornstein--Uhlenbeck diffusion process.
This model has since been considered by several authors, including
\citet{Myers-1981}, \citet{Yashin-1985}, \citet
{Yashin-Vaupel-1986}
and \citet{Aalen-Gjessing-2004}. For given values of the parameters
of the Ornstein--Uhlenbeck process, survival distributions and
hazards are studied. \citet{Myers-1981} focuses on survival
distributions conditioned on initial covariate values;
\citet{Yashin-1985} and \citet{Yashin-Vaupel-1986} use
hazards based
on quadratic functions of Ornstein--Uhlenbeck processes in order to
model heterogeneity among groups and individuals, and to study
the relative hazard functions and survival distributions;
\citet{Aalen-Gjessing-2004} derives
quasi-stationary distributions. Obtaining
such ana\-ly\-ti\-cal results for hazard functions other than
quadratic functions, or for more complex diffusion
processes, is not
feasible.

In our paper, we adopt a Bayesian approach and show how these
models can be efficiently treated by means of Markov chain Monte
Carlo techniques for general choices of diffusion processes and
hazard functions. For instance, by the proposed methods, it is
possible to deal with latent diffusion models which are stochastic
perturbations of common survival models. We also consider the case
of multiple groups of observations, typical of clinical trials, and
we show how to efficiently deal with covariates. We illustrate the
methods via simulation studies and applications to real-world data.

It should be mentioned that other classes of Bayesian nonparametric
and semi-parametric models for survival analysis have been proposed
in the literature. Among the most important, we mention the models based
on \textit{neutral to the right random probabilities}, whose
cumulative hazard rates are processes with independent increments
(see \citet{Doksum-1974} and \citet{Ferguson-1974} for the
definition and
properties of these random measures, and, e.g.,
\citet{Susarla-VanRyzin-1976}, \citet{Kalbfleisch-1978},
\citet{Ferguson-Phadia-1979}, \citet{Hjort-1990} and
\citet{Damien-Walker-2002} for applications in survival analysis),
and all models falling within the framework of \textit{multiplicative
intensity models}, whose hazard rates are mixtures of known kernels
where the mixing measure is a weighted gamma process (see
\citet{Dykstra-Laud-1981}, \citet{Lo-Weng-1989},
\citet{Ishwaran-James-2004} and references therein).

The paper is organized as follows. In Section~\ref{prelim}, we
recall the essentials of diffusion processes and introduce the
model; we also outline how, in the described framework, it is
possible to consider stochastic perturbations of common survival
models. In Section~\ref{McMC-algorithm}, we describe the MCMC scheme
and gives the details of a suitable Hastings-within-Gibbs algorithm,
showing its implementation by means of a toy example.
In Section~\ref{pnc-section}, we present improved versions of the algorithm,
based on reparametrizations of the model. In Section~\ref
{sec:multiple_groups}, we discuss a straightforward
generalization of the framework developed in the previous sections
and deal with the case of multiple groups of observations; this is
also illustrated by application to a data set from a clinical trial, one
that has been considered in a number of papers in the context of
survival analysis, the famous paper by Cox (\citeyear{Cox-1972}) being
among the earliest. In
Section \ref{sec:covariates}, we describe how covariates can be
efficiently included in the proposed models and give an
illustrative application to the lung cancer data set analyzed by
\citet{Muers-Shevlin-Brown-1996}. Finally, in Sections
\ref{sigma-unknown} and \ref{sec:discussion}, we discuss possible
extensions of the models considered.

\section{Latent diffusion models}
\label{prelim}

Let $\Theta$ be a random variable with values in $\mathbb{R}^d$.
Denote by $C([0,\infty),\mathbb{R})$ the space of continuous
functions from $[0,\infty)$ to $\mathbb{R}$ and by $\mathcal{C}$
its cylinder $\sigma$-algebra. Given $\Theta=\theta$, consider the
scalar diffusion process $X=\{X_{t}\dvt t\geq0\}$,  solution of
a \textit{stochastic differential equation} (SDE) of the
form
\begin{eqnarray} \label{SDE}
\mathrm{d}X_t&=&\beta(X_t,\theta)\,\mathrm{d}t+\sigma \,\mathrm{d}B_t,\qquad t\geq0,\nonumber
\\[-8pt]\\[-8pt]
X_0&=&x_0,\nonumber
\end{eqnarray}
driven by the standard scalar Brownian motion $B=\{B_t\dvt t \geq
0\}$. The Brownian motion $B$ and the diffusion process $X$ are
random elements of $(C([0,\infty),\mathbb{R}),\mathcal{C})$. The
diffusion coefficient $\sigma$ is assumed constant and known, for
the moment. The more technically difficult case of unknown $\sigma$
is postponed to Section~\ref{sigma-unknown}. The drift
$\beta(x,\theta)$ is assumed to be
jointly measurable in $x$ and $\theta$, and to
satisfy the regularity conditions (locally Lipschitz, with linear
growth bound) that guarantee the existence of a weakly unique global
solution to (\ref{SDE}). See, for example
\citet{Rogers-Williams-2000b}, Chapter V.24.

Let $\mathbb{W}_{\sigma}$ be the law of $\sigma B$ and, for a given
$\theta$, denote by $\mathbb{P}_{\theta}$ the law of the diffusion
$X$, solution of (\ref{SDE}). By Girsanov's theorem, the
Radon--Nikodym derivative of $\mathbb{P}_{\theta}$ with respect to
$\mathbb{W}_{\sigma}$ is given by
\[
\frac{\mathrm{d}\mathbb{P}_{\theta}}{\mathrm{d}\mathbb{W}_{\sigma}}(x)=\exp\biggl\{\int
_0^{\infty}\frac{\beta(x_t,\theta)}{\sigma^2}\,\mathrm{d}x_t
-\frac{1}{2}\int_0^{\infty}\frac{\beta(x_t,\theta)^2}{\sigma^2}\,
\mathrm{d}t\biggr\},
\]
where $x$ is an element of $(C([0,\infty),\mathbb{R}),\mathcal{C})$.
See, for example, \citet{Rogers-Williams-2000b}, Chapter V.27.

Similarly, for a finite $T$, denote by $C([0,T],\mathbb{R})$ the
space of continuous functions from $[0,T]$ to $\mathbb{R}$ and by
$\mathcal{C}^T$ its cylinder $\sigma$-algebra. Then,
$B_{[0,T]}:=\{B_t\dvt 0\leq t\leq T\}$ and
$X_{[0,T]}=\{X_{t}\dvt 0\leq t\leq T\}$ are random elements of
$(C([0,T],\mathbb{R}),\mathcal{C}^T)$. Let $\mathbb{W}_{T,\sigma}$
be the law of $\sigma B_{[0,T]}$ and, for a given $\theta$, denote
by $\mathbb{P}_{T,\theta}$ the law of $X_{[0,T]}$. Then, by
Girsanov's theorem, the Radon--Nikodym derivative of
$\mathbb{P}_{T,\theta}$ with respect to $\mathbb{W}_{T,\sigma}$ is
given by
%
\begin{equation}
\label{girsanov}
\frac{\mathrm{d}\mathbb{P}_{T,\theta}}{\mathrm{d}\mathbb{W}_{T,\sigma
}}\bigl(x_{[0,T]}\bigr)=\exp\biggl\{\int_0^T\frac{\beta(x_t,\theta)}{\sigma^2}\,\mathrm{d}x_t
-\frac{1}{2}\int_0^T\frac{\beta(x_t,\theta)^2}{\sigma^2}\,\mathrm{d}t\biggr\}
\end{equation}
and, for each $T$, the measures $\mathbb{P}_{T,\theta}$ are
absolutely continuous.

Given the diffusion $X$, let us consider the random distribution
function $F_{X,h}$ on $[0,\infty)$, defined as
%
\begin{equation}
\label{DDRM}
F_{X,h}(t):=1-\exp\biggl\{-\int_0^t h(X_s)\, \mathrm{d}s\biggr\},\qquad t\geq0,
\end{equation}
where $h(\cdot)$ is some suitable non-negative and continuous
function with $\int_0^{\infty} h(X_s)\,\mathrm{d}s=\infty$ almost surely.
The function $h(\cdot)$ plays the role of the hazard function and
$h(X_t)$ is the random hazard rate at time $t$ associated with the
random distribution
$F_{X,h}$.

Two features of the random measure $F_{X,h}$ have to be noted. The
first is
that the hazard inherits the Markov property of the diffusion
process so that the hazard at a future time $t'$ depends only on
the hazard at the present time $t$. Indeed, the Markov property seems
a sensible choice to make at the level of the hazard. The
second is that the cumulative hazard is a process with positively
correlated increments, being the integral of a con\-ti\-nuous
process. The latter feature is natural in many contexts and it
introduces to the model a concern with the stochastic process that
clearly must lie behind the occurrence of events. In words, a high
increment of the cumulative hazard over the time interval $[t,t']$
means that the underlying stochastic process has reached a region of
high risk and this is likely to yield a high increment of the
cumulative hazard over a close (disjoint) time interval. The
strength of this positive correlation, and thus the smoothness of
the cumulative hazard, depends on the choice of the hazard function
$h$ and of the diffusion process~$X$: the rougher the diffusion, the
weaker the correlation, and vice versa; see also the comments in
Section~\ref{sec:discussion}. Note that the property we have just
highlighted differentiates the models we are considering from models
based on neutral to the right random probabilities, whose cumulative
hazards are processes with independent increments and thus have an
erratic behaviour.

Let us now consider a sequence of event times $Y_1,Y_2,\ldots$
which are, conditionally on $F_{X,h}$, independent and identically
distributed (i.i.d.)~with common distribution $F_{X,h}$.
From (\ref{DDRM}), it follows that the distribution of
$Y_1,\ldots,Y_n$, given $X=x$, has density, with respect to the
$n$-dimensional Lebesgue measure $\mathcal{L}^n$, given by
%
\begin{equation}\label{likelyhood}
l(y_1,\ldots,y_n\vert x):=\Biggl[\prod_{j=1}^{n}
h(x_{y_j})\Biggr]\exp\Biggl\{-\sum_{j=1}^{n}\int_0^{y_j} h(x_t)
\,\mathrm{d}t\Biggr\}.
\end{equation}
Censored observations can easily be dealt with in this setting. In
the present paper, we shall restrict our attention to independent
right-censored schemes. If we let $(y_1,\ldots,y_m)$ be the observed event
times and let $(y_{m+1}\!+,\ldots,y_n\!+)$ be the right-censored
event times, then the likelihood becomes
\begin{eqnarray*}
&&l(y_1,\ldots,y_m,y_{m+1}\!+,\ldots,y_n \! +\vert x)
\\
&&\quad=\Biggl[\prod_{j=1}^{m}
h(x_{y_j})\Biggr]\exp\Biggl\{-\sum_{j=1}^{m}\int_0^{y_j} h(x_t)
\,\mathrm{d}t-  \sum_{j=m+1}^{n}\int_0^{y_j +} h(x_t)
\,\mathrm{d}t\Biggr\}.
\end{eqnarray*}

We are thus considering a latent diffusion model for survival
analysis, where the survival times are modeled via an
underlying diffusion process which determines the hazard rate. As
highlighted by \citet{Aalen-Gjessing-2004}, this model can also
be interpreted as a random barrier hitting model. Indeed, the event
occurs when the cumulative hazard strikes a random barrier $R$,
which is exponentially distributed with mean 1 and is
stochastically independent of $X$.

\subsection{Stochastic perturbations of common survival models}
\label{sec:perturbations_common_models}

In the framework we have described, one possibility is to consider
stochastic perturbations of common survival models. Heuristically,
the idea is that if we can express the hazard $r(t)$ of a given
model as a solution of an ordinary differential equation
$\frac{\mathrm{d}r(t)}{\mathrm{d}t}=g(r(t))$ for some suitable function $g$, then
we may be able to use $g$, or some modification of it, to model the
drift of an SDE. Starting from this SDE, we can thus consider a latent
diffusion model whose hazard function is a stochastic perturbation
of $r(t)$.

We shall illustrate this by means of some examples. The simplest case is
offered by the Gompertz model. The Gompertz hazard
$r(t)=\beta\exp\{\alpha t\}$, for $\alpha,\beta > 0$, is a
solution of the\vspace{1pt} ordinary differential equation
$\frac{\mathrm{d}r(t)}{\mathrm{d}t}=g(r(t))=\alpha r(t)$. Consider, thus, the latent
diffusion model based on the SDE having drift $g(X_t)=\theta X_t$
for $\theta>0$,
%
\begin{equation}
\label{eq:SDE_gompertz}
 \mathrm{d}X_t=\theta X_t\,\mathrm{d}t+\sigma\, \mathrm{d}B_t,\qquad t\geq0,\qquad
X_0=x_0>0,
\end{equation}
and with hazard function $h(u)=|u|$. For $\sigma=0$, the SDE
(\ref{eq:SDE_gompertz}) reduces to the ordinary differential equation
written above, for which the Gompertz hazard is a solution, and the
latent diffusion model reduces to the Gompertz model. Hence, the
latent diffusion model based on the SDE~(\ref{eq:SDE_gompertz})
with hazard function $h(u)=|u|$ can be seen as a stochastic
perturbation around a central Gompertz model. This constitutes a
simple example of a latent diffusion model, for which the law of $X_t$,
and thus also the law of the hazard, is known. In the other examples
we shall now give, the SDE cannot be explicitly solved, but the
latent diffusion models based on them can be treated by the
techniques described in the present paper.

Let us consider the Weibull model, whose hazard $r(t)=\alpha\beta
t^{\alpha-1}$ for $\alpha,\beta >  0$ is a non-trivial solution
of the ordinary differential equation
$\frac{\mathrm{d}r(t)}{\mathrm{d}t}=g(r(t))=\gamma r(t)^{(\alpha-2)/(\alpha-1)}$.
Consider, thus, the latent diffusion model based on the SDE
%
\begin{equation}
 \label{eq:SDE_weibull}
 \mathrm{d}X_t=\theta_1 (\operatorname{sign}(X_t)) |X_t|^{\theta_2}\, \mathrm{d}t + \sigma
\,\mathrm{d}B_t,\qquad t\geq
0,\qquad
X_0=x_0>0,
\end{equation}
where
\begin{eqnarray*}
\operatorname{sign}(u)=
\cases{ 1,&\quad \mbox{if} $u>0,$
\cr
-1,&\quad \mbox{if} $u<0,$
\cr
 0,&\quad  \mbox{if} $u=0,$
}
\end{eqnarray*}
and with hazard function $h(u)=|u|$. For $\sigma=0$, the SDE
(\ref{eq:SDE_weibull}) reduces to the ordinary differential equation
written above, for which the Weibull hazard is a solution
($\theta_2$ here plays the role of $(\alpha-2)/(\alpha-1)$). Hence,
the latent diffusion model based on the SDE (\ref{eq:SDE_weibull}),
with hazard function $h(u)=|u|$, can be seen as a stochastic
perturbation around a central Weibull model. For values of
$\theta_2$ in the interval $(0,1)$, which correspond to $\alpha>2$,
the SDE (\ref{eq:SDE_weibull}) has a non-explosive solution. This
solution is weakly unique (see, e.g., \citet{Stroock-Varadhan-2006}).
In Sections \ref{sec:realdata}
and~\ref{sec:appl_covariates}, we shall implement this latent diffusion
model in some illustrative applications to real-world data.

Using the simple idea outlined above, it is possible to develop
other latent diffusion models, such as stochastic perturbations of
log-logistic models and exponential-power models. The log-logistic
hazard ($r(t)=\alpha\beta t^{\alpha-1}/(1+\beta t^{\alpha})$ for
$\alpha,\beta >  0$) and
the exponential-power hazard ($r(t)=\alpha\beta^{\alpha} t^{\alpha
-1}\exp\{(\beta t)^{\alpha}\}$
for $\alpha,\beta >  0$) can, in fact, be written
as solutions of $\frac{\mathrm{d}r(t)}{\mathrm{d}t}=g(r(t))$ for suitable functions
$g$ (when $\alpha < 1$ for the log-logistic
and $\alpha > 1$ for the exponential-power). Let us give a further
example, which generalizes the Pareto model. The Pareto hazard
$r(t)=\alpha/t$, for $\alpha > 0$ and $t \geq \lambda >  0$,
is a solution of the equation
$\frac{\mathrm{d}r(t)}{\mathrm{d}t}=g(r(t))=-\frac{1}{\alpha}[r(t)]^2$. Now, the SDE
having drift $g(X_t)=-\theta X_t^2$, for $\theta>0$,
%
\begin{equation}
\label{SDE_pareto_expl}
 \mathrm{d}X_t=-\theta X_t^2\,\mathrm{d}t+\sigma \,\mathrm{d}B_t,\qquad t\geq\lambda>
0,\qquad X_{\lambda}=x_{\lambda}>0,
\end{equation}
provides a stochastic perturbation around the Pareto hazard, but,
unfortunately, this SDE cannot be used for our purposes since it has
an explosive solution. On the other hand, we can modify
(\ref{SDE_pareto_expl}), for example, by inclusion of $X_t$ in the
diffusion coefficient, in order to obtain another SDE,
%
\begin{equation}
\label{SDE_pareto}
 \mathrm{d}X_t=-\theta X_t^2\,\mathrm{d}t+\sigma X_t \, \mathrm{d}B_t,\qquad t\geq\lambda>
0,\qquad X_{\lambda}=x_{\lambda}>0,
\end{equation}
that also provides a stochastic perturbation around the Pareto
hazard, but has a non-explosive solution. The latter SDE can thus be
transformed into one of constant diffusion coefficient, which can, in
turn, be used in the latent diffusion model. Note that the solution
of (\ref{SDE_pareto}), and that of the corresponding SDE with constant
coefficient, are almost surely positive and so we can take as
hazard function $h(\cdot)$ the identity function, obtaining a
particularly natural perturbation of the Pareto. It is worth
recalling that an SDE with general diffusion coefficient
$\sigma(X_t,\theta)$,
\begin{eqnarray*}
 \mathrm{d}X_t=\beta(X_t,\theta)\,\mathrm{d}t+\sigma(X_t,\theta) \,\mathrm{d}B_t,\qquad
t\geq
0,\qquad X_0=x_0,
\end{eqnarray*}
can, in fact, be transformed into an SDE of unit diffusion coefficient
for the process $Y$, by applying the 1--1 transformation $X_t \to
\eta(X_t;\theta) =: Y_t$, where $\eta(x; \theta) =
\int^{x}\frac{1}{\sigma(z; \theta)} \,\mathrm{d}z$ is any anti-derivative of
$\sigma^{-1}(\cdot; \theta)$ (we are assuming that $\sigma(x,\theta
)$ is
differentiable for any $x\in C([0,\infty),\mathbb{R})$); see, for
example, \citet{beskos-papaspiliopoulos-roberts-fearnhead-2006}. This approach
opens up to a number of possible stochastic perturbations of
commonly used hazards.

\section{Markov chain Monte Carlo methods for latent diffusion
models} \label{McMC-algorithm}

Let $p_{\Theta}(\theta)$ be the prior density, with respect to
$\mathcal{L}^d$, of the $d$-dimensional parameter $\Theta$ which
appears in the drift of the diffusion process $X$, solution of
(\ref{SDE}). Fix a finite time horizon $T$ of interest, with $T\geq
y_{[n]}$, where $y_{[n]}:=\max\{y_1,\ldots,y_n\}$. The choice of $T$
will be discussed in Section \ref{pnc-section}. Then, the joint
posterior distribution of $\Theta$ and $X_{[0,T]}$ has density, with
respect to the product measure
$\mathcal{L}^d\otimes\mathbb{W}_{T,\sigma}$,
given by
%
\begin{equation}
\label{posterior}
\pi\bigl(\theta,x_{[0,T]}\vert y_1,\ldots,y_n\bigr)=C p_{\Theta}(\theta)
g\bigl(x_{[0,T]}\vert \theta\bigr) l\bigl(y_1,\ldots,y_n\vert x_{[0,y_{[n]}]}\bigr),
\end{equation}
where $C$ is a normalizing constant and
$g(x_{[0,T]}|\theta):=\frac{\mathrm{d}\mathbb{P}_{T,\theta}}{\mathrm{d}\mathbb
{W}_{T,\sigma}}(x)$
is given by Girsanov's formula~(\ref{girsanov}).

A Gibbs sampling algorithm for sampling from (\ref{posterior})
alternates between
\begin{enumerate}[1.]
\item simulation of $\Theta$,
conditional on the observations
and the current path of $X_{[0,T]}$;
\item simulation of $X_{[0,T]}$, conditional on the
observations and the current value of $\Theta$.
\end{enumerate}

Note that the parameter $\Theta$ and the observations
$Y_1,\ldots,Y_n$ are conditionally independent, given the
non-observed process $X_{[0,T]}$. In particular, from
(\ref{posterior}), the conditional distribution of $\Theta$ given
$X_{[0,T]}$ has density, with respect to $\mathcal{L}^d$,
proportional to $p_{\Theta}(\theta) g(x_{[0,T]}\vert \theta)$.
The update of the
parameter is particularly straightforward when a conjugate prior
$p_{\Theta}(\theta)$ is
chosen so that it is possible to analytically derive the
conditional distribution of $\Theta$ given $X_{[0,T]}$ and sample
directly from it. The second step is computationally more demanding.
From (\ref{posterior}), the conditional distribution of $X_{[0,T]}$,
given parameter and observations, has density, with respect to
$\mathbb{W}_{T,\sigma}$, proportional to $g(x_{[0,T]}\vert \theta)
l(y_1,\ldots,y_n\vert x)$ and cannot be sampled directly. An appropriate
Metropolis--Hastings step is thus required.

Implementation of the algorithm will necessary involve a
discretization of the diffusion sample path. When the SDE cannot be
solved, it is possible to use \textit{Euler--Maruyama approximation};
see, for example, Chapter 9 in \citet{Kloeden-Platen-1992}.
Alternatively, it may be possible to simulate the diffusion path by
means of the exact algorithm described in
\citet{beskos-papaspiliopoulos-roberts-fearnhead-2006}, thus
avoiding approximation errors.

\subsection{Hastings-within-Gibbs algorithm for a latent diffusion model}\label{McMC-section}

We now give the details of the Hastings-within-Gibbs algorithm for
latent diffusion models.

Just as an example, consider a latent diffusion model with base
diffusion which is  solution of the SDE
%
\begin{equation}
\label{SDE-model-general}
 \mathrm{d}X_t=\theta^{\mathsf{ T}} f(X_t)\,\mathrm{d}t+\sigma\, \mathrm{d}B_t,\qquad
t\geq0,\qquad X_0=x_0,
\end{equation}
with $\theta^{\mathsf{T}}=(\theta_1, \ldots,\theta_d)$ and
$f(x)^{\mathsf{ T}}=(f_1(x),\ldots,f_d(x))$, where $f_i(x)$
is some real-valued function for $i=1,\ldots,d$. Let the drift
$\theta^{\mathsf{ T}} f(x)$ be such that the regularity
conditions mentioned in Section \ref{prelim} are satisfied. Let the
prior for $\Theta=(\Theta_1,\ldots,\Theta_d)$ be multivariate
Gaussian with mean vector and variance matrix, respectively, given by
\begin{eqnarray*}
\mu=\left[
\matrix{\mu_1\cr
\mu_2\cr
\vdots\cr
\mu_d
}
\right]
\quad\mbox{and}\quad\Sigma=\left[
\matrix{
\lambda_{11}&\lambda_{12}&\cdots&
\lambda_{1d}\cr
\lambda_{12}&\lambda_{22}&\cdots&
\lambda_{2d}\cr
\vdots&\vdots&\ddots&\vdots\cr
\lambda_{1d}&\lambda
_{2d}&\cdots&
\lambda_{dd}
}
\right]^{-1}.
\end{eqnarray*}
Then, the distribution of $\Theta$,
given the diffusion $X_{[0,T]}=x_{[0,T]}$, is still Gaussian, with
mean and covariance matrix, respectively, given by
%
\begin{equation}\label{mean-var-full-cond1}
\mu_x=\Sigma_x\left[
\matrix{S_1\cr
S_2\cr
\vdots\cr
S_d
}
\right]
\quad\mbox{and}\quad\Sigma_x=\left[
\matrix{ L_{11}&L_{12}&\cdots&
L_{1d}\cr
L_{12}&L_{22}&\cdots&
L_{2d}\cr
\vdots&\vdots&\ddots&\vdots\cr
L_{1d}&L_{2d}&\cdots&
L_{dd}
}
\right]^{-1},
\end{equation}
where, for $i=1,\ldots,d$ and $j=1,\ldots,d$,
\begin{eqnarray*}
S_i &:=& \frac{1}{\sigma^2}\int_0^Tf_i(x_t)\,\mathrm{d}x_t+\sum_{j=1}^{d}\lambda
_{ij}\mu_j,\qquad
L_{ij}:=\frac{1}{\sigma^2}\int_0^T f_i(x_t)f_j(x_t)\,\mathrm{d}t+\lambda_{ij}.
\end{eqnarray*}
The update of $\Theta$ can thus be performed by sampling directly
from this conditional distribution.

The update of the diffusion $X_{[0,T]}$ is less straightforward and
requires an appropriate Metropolis--Hastings step. It is possible, for
example, to carry out an independence sampler with proposal
distribution given by a Brownian motion starting at $x_0$. To
improve the acceptance rate of the move that updates the diffusion
path, we apply the following updating strategy. Let
$0=t_1<\cdots<t_m=T$. Instead of proposing a new diffusion path on
the whole interval $[0,T]$, we propose to change the trajectory only
on a subinterval $[t_i,t_{i+2}]$, keeping the rest of the
diffusion fixed. To ensure continuity of the diffusion path, the proposal
distribution for the new trajectory on the subinterval
$[t_i,t_{i+2}]$ is a Brownian bridge
$\mathit{BB}_{[t_i,t_{i+2}]}(x_{t_{i}},x_{t_{i+2}})=\{\mathit{B B}_t(x_{t_{i}},x_{t_{i+2}})\dvt t_i\leq
t\leq t_{i+2}\}$, having
as starting and ending points, respectively, the values $X_{t_i}=x_{t_i}$
and $X_{t_{i+2}}=x_{t_{i+2}}$ of the current diffusion. The proposed
diffusion path $x^*_{[0,T]}$ is then given by $\{x^*_t=1(t
\notin[t_i,t_{i+2}])x_t+ 1(t \in
[t_i,t_{i+2}])bb_t(x_{t_{i}},x_{t_{i+2}})\dvt t \in[0,T]\}$, where
$\mathit{bb}_t(x_{t_{i}},x_{t_{i+2}})$ is the realization of the Brownian
bridge $\mathit{BB}_{[t_i,t_{i+2}]}(x_{t_{i}},x_{t_{i+2}})$. This move is
accepted with probability
%
\begin{equation}\label{acptcoeff}
1\wedge
\frac{g(\mathit{bb}_{[t_i,t_{i+2}]}(x_{t_{i}},x_{t_{i+2}})|\theta
)}{g(x_{[t_i,t_{i+2}]}|\theta)}
\frac{l(y_1,\ldots,y_n|x_{[0,y_{[n]}]}^*)}
{l(y_1,\ldots,y_n|x_{[0,y_{[n]}]})},
\end{equation}
where $g(x_{[t_i,t_{i+2}]}\vert \theta)$ is given by Girsanov's formula
restricted to the interval $[t_i,t_{i+2}]$, that is,
\begin{eqnarray*}
g\bigl(x_{[t_i,t_{i+2}]}\vert \theta\bigr)=\exp\biggl\{\int_{t_i}^{t_{i+2}}\frac{\theta
^{\mathsf{
T}} f(X_t)}{\sigma^2}\,\mathrm{d}x_t-
\frac{1}{2}\int_{t_i}^{t_{i+2}}\frac{(\theta^{\mathsf{
T}} f(X_t))^2}{\sigma^2}\,\mathrm{d}t\biggr\}.
\end{eqnarray*}
The procedure is iterated for $i=1,\ldots,m-3$. Note that the
different blocks $[t_i,t_{i+2}]$ overlap so that there are no time
instants where the diffusion is kept fixed. For the same reason, the
last block $[t_{m-2},T]$ is updated by means of a Brownian motion
$B_{[t_{m-2},T]}(x_{t_{m-2}})$ starting at
$X_{t_{m-2}}=x_{t_{m-2}}$ so that the value of the diffusion at $T$
may vary. The acceptance coefficient of the move that updates the
last block is the same as in (\ref{acptcoeff}), with
$[t_i,t_{i+2}]=[t_{m-2},T]$ and $b_{[t_{m-2},T]}(x_{t_{m-2}})$ in
place of $\mathit{bb}_{[t_i,t_{i+2}]}(x_{t_{i}},x_{t_{i+2}})$, where
$b_{[t_{m-2},T]}(x_{t_{m-2}})$ is the realization of the Brownian
motion $B_{[t_{m-2},T]}(x_{t_{m-2}})$.

This idea of updating smaller intervals at a time has been used in
\citet{Shephard-Pitt-1997}
for the simulation
of non-Gaussian time series models and later applied for the
simulation of discretely observed diffusions, for example, by
\citet{elerian-chib-shephard-2001}.

In Section \ref{simulation-p-e}, we shall illustrate the
implementation of this algorithm by means of a toy example. Note
that in this section and in the following, we are considering base
diffusions having drift linear in the parameter $\theta$ simply for
purposes of
exposition.

\subsection{Implementation of the algorithm: A toy example} \label
{simulation-p-e}

\begin{figure}[t]

\includegraphics{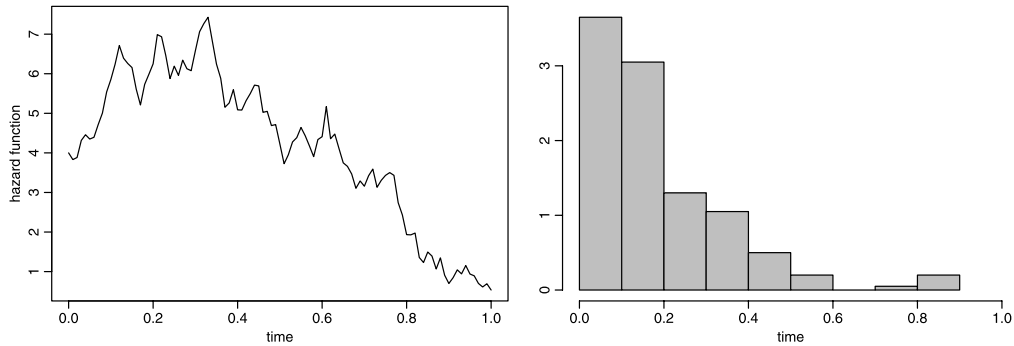}

\caption{Left: hazard function $x^2$. Right: histogram of
data sampled from $F_{x,x^2}$ with censoring at $C=0.9$.}
\label{fig:descrittiva}
\end{figure}

We show here the implementation of the algorithm described in
Section \ref{McMC-section}, by means of a toy example. Consider the
model based on the diffusion process satisfying the SDE
%
\begin{equation}
\label{SDE-model}
 \mathrm{d}X_t=\theta_1 \sin(X_t)\,\mathrm{d}t+\theta_2 \,\mathrm{d}t+ \mathrm{d}B_t,\qquad t\geq
0,\qquad X_0=2,
\end{equation}
with hazard function $h(u)=u^2$. We simulate observations from this
model for values of the parameters $\theta_1=-1.4$ and
$\theta_2=-1$, and censoring time $C=0.9$. In particular, we
sample one realization $x$ of the diffusion process satisfying
(\ref{SDE-model}), with $\theta_1=-1.4$ and $\theta_2=-1$. We
then simulate 200 i.i.d. observations from the corresponding
distribution $F_{x,h}=1-\exp\{-\int_0^t (x_s)^2
\,\mathrm{d}s\}$ and censor the observations at a common cut-off
$C=0.9$. The diffusion is sampled at intervals of length $0.01$,
using Euler--Maruyama approximation. Figure \ref{fig:descrittiva}
shows the corresponding hazards (the squared diffusion) and a
histogram of sampled data. The hazard function has a typical shape,
first (mainly) increasing and then (mainly) decreasing.

We choose as time horizon of interest $T=1$. We then run the
Hastings-within-Gibbs algorithm under the following specifications.
The prior for $(\theta_1,\theta_2)$ is Gaussian, as in Section
\ref{McMC-section}, with $\mu_1=-1.4$, $\mu_2=-1$,
$\lambda_{11}=\lambda_{22}=1/5$ and $\lambda_{12}=0$. The starting
values of the parameters are $\theta_1=\theta_2=0$ and the starting
diffusion is a Brownian motion, starting at $x_0=2$. The diffusion
path is updated on subintervals of length 0.2 at a time.
The algorithm is run for 200\,000 iterations and the
first 2000 are discarded as burn-in.

\begin{figure}[t]

\includegraphics{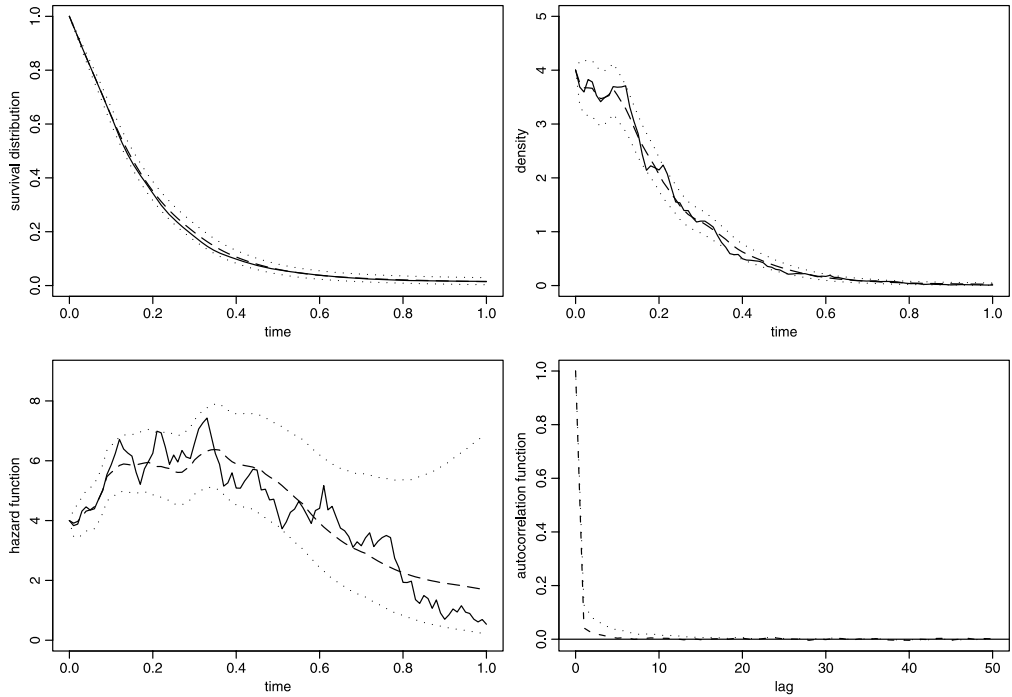}

\caption{Upper-left:
true survival distribution $1-F_{x,x^2}$ (solid), together with its
posterior mean (dashed) and pointwise approximate $90\%$ highest
posterior bands (dotted). Upper-right: true density (solid), together
with its posterior mean (dashed) and pointwise approximate $90\%$
highest posterior bands (dotted). Lower-left: true hazard function
$x^2$ (solid), together with its posterior mean (dashed) and
pointwise approximate $90\%$ highest posterior bands (dotted).
Lower-right: autocorrelation functions for $\theta_1$ series
(dotted) and $\theta_2$ series (dashed).}
 \label{fig:risultati}
\end{figure}

Figure \ref{fig:risultati} shows the estimates of survival
distribution, density and hazard function, based on the MCMC
output, together with pointwise approximate $90\%$ highest
posterior bands. The true survival distribution and hazard function
are also displayed to demonstrate the good fit of the MCMC estimates.
Figure \ref{fig:risultati} also shows autocorrelation functions for
$\theta_1$ and $\theta_2$ series.


\section{Reparametrizations of the latent diffusion models} \label
{pnc-section}

\begin{figure}[t]

\includegraphics{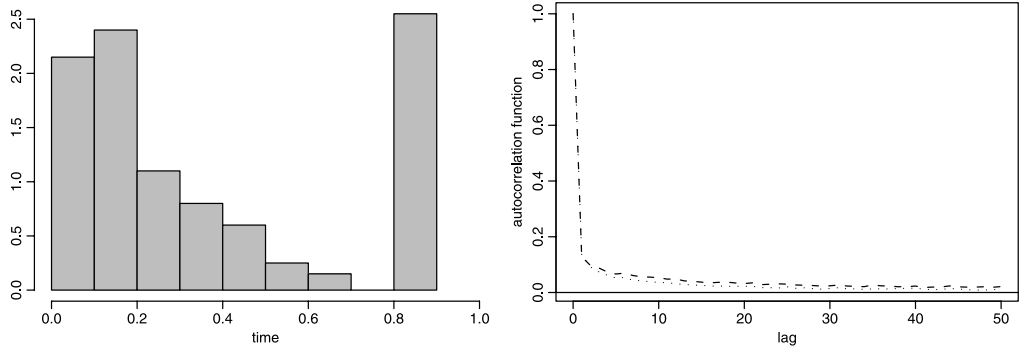}

\caption{Left: histogram
of data sampled from $F_{x' ,x^{\prime2}}$ with censoring at $C=0.9$. Right:
autocorrelation functions for $\theta_1$ series (dotted) and
$\theta_2$ series (dashed).}
\label{fig:descrittiva-c-pnc}
\end{figure}

The MCMC algorithm described in the previous sections might have
poor mixing properties when we consider a finite-time horizon $T$
significantly larger than the maximum of the data. This problem is
evident in Figure \ref{fig:descrittiva-c-pnc}. This figure shows
the histogram of 200 i.i.d. observations from the distribution
$F_{x' ,h}$, where $x'$ is a new realization of the diffusion
process satisfying the same SDE used in Section
\ref{simulation-p-e}; also, the hazard function $h$ and the
censoring time $C$ are the same. In this simulation, we have fixed a
longer time horizon $T=1.8$ and have then run the algorithm
under the same specifications of Section \ref{simulation-p-e}.
Figure \ref{fig:descrittiva-c-pnc} displays autocorrelation
functions for $\theta_1$ and $\theta_2$ series, which are not
exponentially decreasing. With the same data set, but choosing a
shorter time horizon (such as $T=1$, as in the previous section),
the algorithm does not exhibit strong serial correlation in the
draws of $\theta_1$ and $\theta_2$. The worsening of the mixing
properties of the algorithm when $T$ becomes significantly larger
than the maximum of the data was also observed for the data set
simulated in Section \ref{simulation-p-e}.

To avoid this problem, we propose a modification of the algorithm
which has good mixing properties, regardless of the choice of time
horizon, and is, in fact, completely robust with respect to $T$. The
algorithm is based on a simple reparametrization of the model.
Indeed, the performance of MCMC methods, particularly when using
Gibbs samplers, depends crucially on the parametrization of the
unknown quantities in the hierarchical structure. The issue of
reparametrization of the posterior distributions in order to improve
convergence properties of the algorithms has received much
attention. See, for example, \citet{Hills-Smith-1992},
\citet{Gelfand-Sahu-Carlin-1995}, \citet{Gelfand-Sahu-Carlin-1996}
and
\citeauthor{Papaspiliopoulos-Roberts-Skold-2003}
(\citeyear{Papaspiliopoulos-Roberts-Skold-2003,papaspiliopoulos-roberts-skold-2007}).

Instead of using the natural parametrization of the model in terms
of $(\Theta,X)$, the so-called \textit{centered parametrization}, we
parametrize it in terms of $(\Theta,\widetilde{X})$, where
\begin{eqnarray*}
\widetilde{X}_t=1\bigl(t\leq
y_{[n]}\bigr)X_t+1\bigl(t>
y_{[n]}\bigr)\bigl[B_t-B_{ y_{[n]}}\bigr],\qquad t\geq0.
\end{eqnarray*}
In the terminology used by
\citet{Papaspiliopoulos-Roberts-Skold-2003}, this is called a
\textit{partially non-centered parametrization}, the fully
\textit{non-centered parametrization} being, in this case,
$(\Theta,B)$.
The diffusion $X$ can then be reconstructed as
a function of $\Theta$, $\widetilde{X}$ and $y_1,\ldots,y_n$, by
\begin{eqnarray*}
\cases{
 X_t=\widetilde{X}_t, &\quad $0\leq t\leq y_{[n]},$
 \cr
\mathrm{d}X_t=\beta(X_t,\Theta)\,\mathrm{d}t+\sigma \,\mathrm{d}\widetilde{X}_t,& \quad $t\geq
y_{[n]}.$
}
\end{eqnarray*}
The joint posterior distribution of $\Theta$ and $\widetilde{X}$ has
density, with respect to the product measure
$\mathcal{L}^d\otimes\mathbb{W}_{\sigma}$,
given by
%
\begin{equation}
\label{posterior-pnc}
\pi(\theta,\tilde{x}\vert y_1,\ldots,y_n)=C
p_{\Theta}(\theta)
g\bigl({x}_{[0,y_{[n]}]}\vert \theta\bigr)
l\bigl(y_1,\ldots,y_n\vert {x}_{[0,y_{[n]}]}\bigr),
\end{equation}
where $x_{[0,y_{[n]}]}\equiv\tilde{x}_{[0,y_{[n]}]}$, $C$ is a
normalizing constant and $g({x}_{[0,y_{[n]}]}\vert \theta)=
\frac{\mathrm{d}\mathbb{P}_{y_{[n]},\theta}}{\mathrm{d}\mathbb{W}_{y_{[n]},\sigma
}}(x_{[0,y_{[n]}]})$
is given by Girsanov's formula (\ref{girsanov}). Note, in particular,
that (\ref{posterior-pnc}) characterizes the posterior distribution
of $\widetilde{X}$, and thus the posterior distribution of the
diffusion $X$, over the whole positive half-line. It thus also
highlights that $X_{[0,y_{[n]}]}$ acts as sufficient statistics.

It is possible to simulate from (\ref{posterior-pnc}) by means of a
Gibbs sampler quite similar to the one described in Section
\ref{McMC-section}. However, the algorithm is now completely robust
to the choice of $T$ since the update of the parameter $\Theta$,
conditionally on $\widetilde{X}$, only involves
$\widetilde{X}_{[0,y_{[n]}]}$. In the first step, in fact, we now
simulate $\Theta$ conditionally on $\widetilde{X}_{[0,y_{[n]}]}$. In
the second step, we simulate $\widetilde{X}$ over the time interval
of interest, $[0,T]$, conditionally on $\Theta$ and the observations.
In this case, we use a proposal distribution which is a Brownian
motion starting at $x_0$ over the time interval $[0,y_{[n]}]$ and
a Brownian motion starting at $0$ over the time interval
$[y_{[n]},T]$. On $[0,y_{[n]}]$, we again follow the updating
strategy with the overlapping Brownian bridges that was described in
Section \ref{McMC-section}. When reconstructing the diffusion
${X}_{[0,T]}$ from $\Theta$ and $\widetilde{X}_{[0,T]}$, we are
careful to preserve the continuity of the diffusion path at time
$y_{[n]}$. Details are omitted.

\begin{figure}[t]

\includegraphics{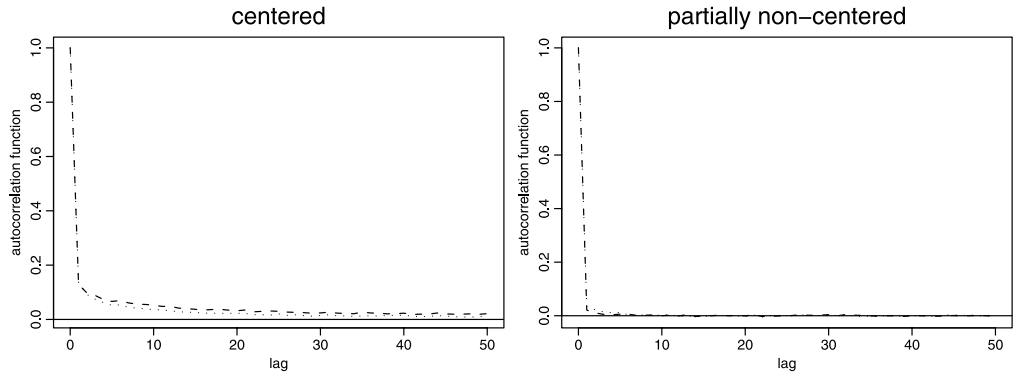}

\caption{Autocorrelation functions
for $\theta_1$ series (dotted) and $\theta_2$ series (dashed),
obtained with the algorithm based on the centered parametrization
(left) and with the algorithm based on the partially non-centered
parametrization (right).} \label{fig:risultati-c-pnc-1}
\end{figure}

\begin{figure}[t!]

\includegraphics{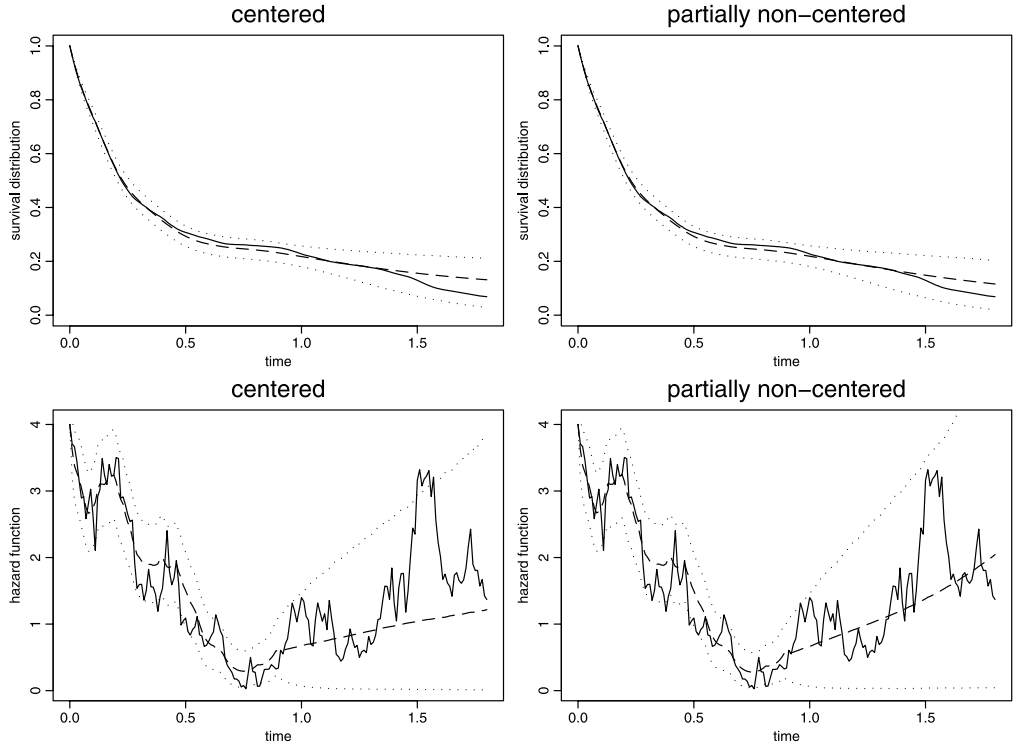}

\caption{Top: true survival distribution distribution
$1-F_{x' ,x'^2}$ (solid), together with its posterior mean
(dashed) and pointwise approximate $90\%$ highest posterior bands
(dotted), obtained with the algorithm based on the centered
parametrization (left) and with the algorithm based on the partially
non-centered parametrization (right). Bottom: true hazard function
$x'^2$ (solid), together with its posterior mean (dashed) and
pointwise approximate $90\%$ highest posterior bands (dotted),
obtained with the algorithm based on the centered parametrization
(left) and with the algorithm based on the partially non-centered
parametrization (right).}
\label{fig:risultati-c-pnc-2}
\end{figure}

Figures \ref{fig:risultati-c-pnc-1} and \ref{fig:risultati-c-pnc-2}
compare mixing and MCMC estimates obtained with the algorithms
based on the centered parametrization and on the partially
non-centered parametrization for the data set corresponding to
Figure \ref{fig:descrittiva-c-pnc}. The specifications of the two
algorithms are as in Section \ref{simulation-p-e}.
Note that the hazard function is bathtub shaped. Hazard functions
with such a shape are quite common in survival analysis (think, for
instance, of human mortality).

%

As we shall see in Section \ref{sec:covariates}, another
reparametrization of the model, one that turns out to be useful in
the presence of covariates, is the fully \textit{non-centered
parametrization} in terms of $(\Theta, B)$. The diffusion $X$ can
be reconstructed as a function of $\Theta$ and $B$, simply by
the SDE
\begin{eqnarray*}
\mathrm{d}X_t=\beta(X_t,\Theta)\,\mathrm{d}t+\sigma \,\mathrm{d}B_t,\qquad t\geq0,\qquad
X_0=x_0.
\end{eqnarray*}
The joint posterior distribution of $\Theta$ and $B$ has density,
with respect to the product measure
$\mathcal{L}^d\otimes\mathbb{W}_{\sigma}$,
given by
%
\begin{equation}
\label{posterior-nc}
\pi(\theta,b\vert y_1,\ldots,y_n)=C
p_{\Theta}(\theta)
l\bigl(y_1,\ldots,y_n\vert \theta,b_{[0,y_{[n]}]}\bigr),
\end{equation}
where $C$ is a normalizing constant and
$l(y_1,\ldots,y_n\vert \theta,b_{[0,y_{[n]}]})=l(y_1,\ldots,y_n\vert x_{[0,y_{[n]}]})$
is as in (\ref{likelyhood}). Note that, similarly to what has been
observed for the partially non-centered parametrization,
(\ref{posterior-nc}) also characterizes the posterior distribution of the
diffusion $X$ over the whole positive half-line. Moreover,
 the Gibbs sampler that simulates from (\ref{posterior-nc})
is also completely robust with respect to the choice of the time horizon
$T$. In the first step, we simulate $\Theta$ conditionally on
$B_{[0,y_{[n]}]}$ and the observations. Note, in particular, that the
conditional distribution of $\Theta$, given $B_{[0,T]}$ and the
observations, now has density, with respect to $\mathcal{L}^d$,
proportional to $p_{\Theta}(\theta)
l(y_1,\ldots,y_n\vert \theta,b_{[0,y_{[n]}]})$. In the second
step, we simulate $B$ over the time interval of interest, $[0,T]$,
conditionally on $\Theta$ and the observations. For proposal
distribution, we use a Brownian motion starting at $0$ and we employ
the updating strategy based on overlapping
Brownian bridges. In this case, when updating the Brownian motion
path $b$ over the subinterval $[t_i,t_{i+2}]$, we need to
reconstruct the corresponding diffusion path $x$ over the
subinterval $[t_i,T]$ in order to preserve the continuity of the
diffusion path at time $t_{i+2}$. Details are omitted.

\section{Latent diffusion models for multiple groups of observations}
\label{sec:multiple_groups}

We now discuss a straightforward generalization of the framework
developed in the previous sections and deal with the case of
multiple groups of observations, where the observations within each
group are taken under homogeneous conditions. Consider, for example,
the case in which different treatments are being administered to
different groups of patients in a clinical trial.

Given $\Theta=\theta$, let $X^{[1]}, \ldots,X^{[q]}$ be $q$
stochastically independent diffusion processes satisfying
(\ref{SDE}) and $F_{X^{[1]},h}, \ldots,F_{X^{[q]},h}$ the relative
random distributions, as in (\ref{DDRM}). Now, consider $q$ sequences
of observations $(Y_n^{[1]})_n, \ldots, (Y_n^{[q]})_n$ such that the
random variables\vspace{1pt} in $((Y_n^{[1]})_n, \ldots,(Y_n^{[q]})_n)$
are conditionally independent, given $F_{X^{[1]},h},
\ldots,F_{X^{[q]},h}$, and the random variables in $(Y_n^{[k]})_n$
have common distribution $F_{X^{[k]},h}$ for $k=1,\ldots,q$.

The joint distribution of
$Y_1^{[1]},\ldots,Y_{n_1}^{[1]},\ldots,$
$Y_1^{[q]},\ldots,Y_{n_q}^{[q]}$, given $X^{[1]}=x^{[1]}, \ldots,$
$X^{[q]}=x^{[q]}$, has density, with respect to $\mathcal{L}^n$
(where $n=n_1+\cdots+n_q$), given by
\begin{eqnarray*}
l\bigl(y_1^{[1]},\ldots,y_{n_1}^{[1]};\ldots;y_1^{[q]},\ldots
,y_{n_q}^{[q]}\big\vert x_{[0,y_{[n_{ 1}]}]}^{[1]},
\ldots, x_{[0,y_{[n_{ q}]}]}^{[q]}\bigr)=
\prod_{k=1}^{q}l\bigl(y_1^{[k]},\ldots,y_{n_k}^{[k]}\big\vert x_{[0,y_{[n_{ k}]}]}^{[k]}\bigr),
\end{eqnarray*}
where $y_{[n_{ k}]}:=\max\{y_1^{[k]},\ldots,y_{n_k}^{[k]}\}$ and
$l(y_1^{[k]},\ldots,y_{n_k}^{[k]}\vert x_{[0,y_{[n_{ k}]}]}^{[k]})$
is as in (\ref{likelyhood}). Using the partially non-centered
parametrization described in Section \ref{pnc-section}, the joint
posterior distribution of $\Theta$ and $\widetilde{X}^{[1]},
\ldots,\widetilde{X}^{[q]}$ has density, with respect to the product
measure $\mathcal{L}^d\otimes\mathbb{W}_{\sigma}^q$, given by
\begin{eqnarray}
\label{posterior-p-e}
&&\pi\bigl(\theta,\tilde{x}^{[1]},
\ldots,\tilde{x}^{[q]}\big\vert y_1^{[1]},\ldots,y_{n_1}^{[1]};\ldots
;y_1^{[q]},\ldots,y_{n_q}^{[q]}\nonumber
\bigr)
\\[-8pt]\\[-8pt]
&&\quad =C p_{\Theta}(\theta)
\Biggl[\prod_{k=1}^{q}g\bigl(x_{[0,y_{[n_{ k}]}]}^{[k]}\vert \theta\bigr)
l\bigl(y_1^{[k]},\ldots,y_{n_k}^{[k]}\big\vert x_{[0,y_{[n_{ k}]}]}^{[k]}\bigr)\Biggr],\nonumber
\end{eqnarray}
where $C$ is a normalizing constant and
$g(x_{[0,y_{[n_{ k}]}]}^{[k]}\vert \theta)=\frac{\mathrm{d}\mathbb{P}_{y_{[n_{
k}]},\theta}}{\mathrm{d}\mathbb{W}_{y_{[n_{ k}]},\sigma}}(x_{[0,y_{[n_{
k}]}]}^{[k]})$
is given by\vspace{-2pt} Girsanov's formula (\ref{girsanov}).

The contributions of the $q$ groups of observations factorize in
(\ref{posterior-p-e}) and a simple modification of the MCMC
algorithm presented in the previous sections may be used to deal
with this case. Let $T_1,\ldots,T_q$ be the time horizons of
interest for the $q$ groups, with $T_k\geq y_{[n_{ k}]}$ for
$k=1,\ldots,q$. The Hastings-within-Gibbs algorithm for sampling
from (\ref{posterior-p-e}) alternates between
\begin{enumerate}[1.]
\item simulation of $\Theta$, conditional on the current paths of
$\widetilde{X}_{[0,y_{[n_{ 1}]}]}^{[1]},
\ldots,$ $\widetilde{X}_{[0,y_{[n_{ q}]}]}^{[q]}$;
\item for each
$k$ in $ \{1,\ldots,q\}$, simulation of
$\widetilde{X}_{[0,T_k]}^{[k]}$, conditional on the observations
$Y_1^{[k]},\ldots,Y_{n_k}^{[k]}$ and the current value of $\Theta$.
\end{enumerate}

Consider, for example, a latent diffusion model with $q$
stochastically independent diffusion processes, $X^{[1]},
\ldots,X^{[q]}$, satisfying the SDE (\ref{SDE-model-general}).
Choose the same multivariate Gaussian prior for $\Theta$ that was used
in Section \ref{McMC-section}. Then, the distribution of
$\Theta$, given
$\widetilde{X}_{[0,y_{[n_{ 1}]}]}^{[1]}=x_{[0,y_{[n_{ 1}]}]}^{[1]},
\ldots,$
$\widetilde{X}_{[0,y_{[n_{ q}]}]}^{[q]}=x_{[0,y_{[n_{ q}]}]}^{[q]}$,
is still Gaussian, with mean vector and covariance matrix as in
(\ref{mean-var-full-cond1}), but with
\begin{eqnarray*}
S_i &:=& \frac{1}{\sigma^2}\Biggl[\sum_{k=1}^q\int_0^{y_{[n_{
k}]}}f_i\bigl(x^{[k]}_t\bigr)\,\mathrm{d}x^{[k]}_t\Biggr]+\sum_{j=1}^{d}\lambda_{ij}\mu_j,
\\
L_{ij}&:=& \frac{1}{\sigma^2}\Biggl[\sum_{k=1}^q\int_0^{y_{[n_{
k}]}}f_i\bigl(x^{[k]}_t\bigr)f_j\bigl(x^{[k]}_t\bigr)\,\mathrm{d}t\Biggr]+\lambda_{ij}
\end{eqnarray*}
for $i=1,\ldots,d$, $j=1,\ldots,d$.
The
update of the parameter $\Theta$ can thus be performed by sampling
directly from this conditional distribution. The second step may be
carried out by $q$ repetitions of the updating mechanism
described in Sections~\ref{McMC-section} and \ref{pnc-section}.

Note that we are here considering a simple hierarchical structure,
where inference on the separate groups is linked only at the
level of the finite-dimensional parameter $\Theta$. For some
applications, this might allow too little borrowing of strength for
inference across groups of patients. In Section \ref{sec:covariates},
we shall instead describe a more complex hierarchical structure,
suitable in the presence of covariates and allowing for a much
stronger borrowing of strength for inference across individuals.

\subsection{An illustrative application to a real data set with
multiple groups of observations}
\label{sec:realdata}

In this section, we show the implementation of the latent diffusion
model for multiple groups of observations via an illustrative
application to a small data set from a clinical trial, one that has been
considered in a number of papers in the context of survival
analysis, among them \citet{Gehan-1965}, \citet{Cox-1972},
\citet{Wei-1984} and \citet{Xu-OQuigley-2000} in the non-Bayesian
literature and \citet{Kalbfleisch-1978},
\citet{Laud-Damien-Smith-1998} and \citet{Damien-Walker-2002} in the
Bayesian literature.
In the trial, reported by
\citet{Freireich-1963}, 6-mercaptopurine (6-MP) was compared to a
placebo in the maintenance of remission in acute leukemia. The
following lengths of remission in weeks were recorded for 42
patients, half of which treated with the 6-MP drug and half with the
placebo (a $+$ sign indicates a censored observation):
\begin{enumerate}[placebo:]
\item[6-MP:] 6, 6, 6, 6$+$, 7, 9$+$, 10, 10$+$, 11$+$, 13, 16, 17$+$, 19$+$, 20$+$,
22, 23, 25$+$, 32$+$, 32$+$, 34$+$, 35$+$,
\item[placebo:]
$1,1,2,2,3,4,4,5,5,8,8,8,8,11,11,12,12,15,17,22,23.$
\end{enumerate}

We thus consider a model for two groups of observations, namely the
6-MP drug group and the placebo group. As latent diffusion model, we
shall use the stochastic perturbation around the Weibull described
in Section \ref{sec:perturbations_common_models}. Recall that this
model has base diffusion sa\-ti\-sfying the SDE
\begin{eqnarray*}
\mathrm{d}X_t=\theta_1 (\operatorname{sign}(X_t)) |X_t|^{\theta_2} \,\mathrm{d}t + \sigma \,\mathrm{d}B_t,\qquad
t\geq0,\qquad
X_0=x_0>0,
\end{eqnarray*}
and hazard function $h(u)=|u|$.

We express the data as fractions of one year and choose as time
horizons of interest $T_1=T_2=0.75$, corresponding to 9 months (39
weeks). We take $\Theta_1$ and $\Theta_2$ to be a priori independent, with
a Gaussian prior distribution for $\Theta_1$, mean $\mu=0$,
variance $1/\lambda=5$, and a uniform prior over $[0,1]$ for
$\Theta_2$. Moreover, we set $x_0=0.8$ and $\sigma=8$. We then run
the Hastings-within-Gibbs algorithm based on the partially
non-centered parametrization. The update of $\Theta_1$ is performed
by sampling directly from the conditional distribution $\Theta_1$,
given
$\Theta_2,\widetilde{X}_{[0,y_{[n_1]}]}^{[1]},\widetilde
{X}_{[0,y_{[n_2]}]}^{[2]}$,\vspace{-2pt}
which is still Gaussian with mean $\frac{S+\lambda\mu}{L+\lambda}$
and variance $\frac{1}{L+\lambda}$, where
\begin{eqnarray*}
S:=\frac{1}{\sigma^2}\Biggl[\sum_{j=1}^2\int_0^{y_{[n_j]}}\bigl(\bigl(\operatorname{sign}\bigl(x^{[j]}_t\bigr)\bigr)
\bigl|x^{[j]}_t\bigr|^{\theta_2}\bigr)\,\mathrm{d}x^{[j]}_t\Biggr]\quad\mbox{and}\quad
L:=\frac{1}{\sigma^2}\Biggl[\sum_{j=1}^2\int
_0^{y_{[n_j]}}\bigl(\bigl|x^{[j]}_t\bigr|^{\theta_2}\bigr)^2\,\mathrm{d}t\Biggr].
\end{eqnarray*}
For the update of $\Theta_2$, we use an independence sampler with a
Beta proposal distribution, with parameters $(1/2,1/2)$.
The update of $\widetilde{X}^{[1]}$ and
$\widetilde{X}^{[2]}$ is carried out as described in the previous
sections.
The algorithm is run for 200\,000
iterations and the
first 2000 are discarded as burn-in.

\begin{figure}[t]

\includegraphics{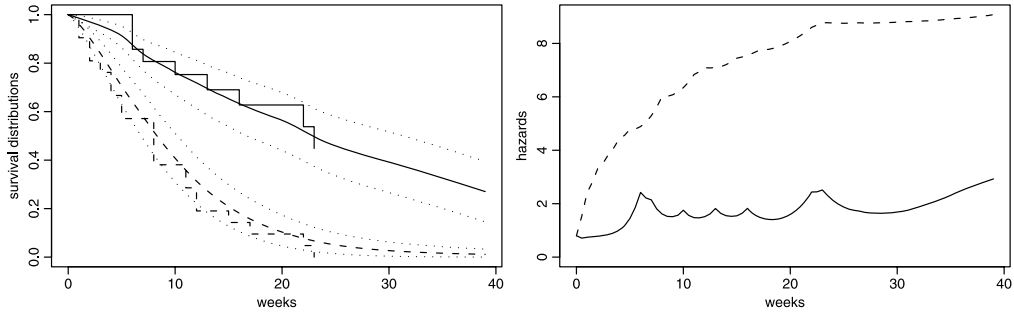}

\caption{Left:
posterior mean survival distributions and pointwise approximate
$90\%$ highest posterior bands for the group of patients treated
with 6-MP drug (solid) and for the group of patients treated with
the placebo (dashed), together with corresponding Kaplan--Meier
curves. Right: posterior mean hazards for the group of patients
treated with 6-MP drug (solid) and for the group of patients treated
with the placebo (dashed).} \label{fig:surv-frereich} \vspace{2mm}
\end{figure}

Figure \ref{fig:surv-frereich} displays the MCMC estimates of the survival
distributions of the two groups, 6-MP drug and placebo, together
with the relative Kaplan--Meier curves. Note that the MCMC estimates
of the two survival distributions are closer to one another than the
two Kaplan--Meier curves, thus indicating borrowing of strength for
inference among the two groups. Hence, the latent diffusion model,
which gains much flexibility over a fully parametric model by
introducing randomness around it, does not suffer from the opposite
problem of being too data-driven.
Figure~\ref{fig:surv-frereich} also
displays the MCMC estimates of the hazards of the two groups.

We could now verify the efficacy of 6-MP drug treatment as proposed
in \citet{Damien-Walker-2002}. In particular, under the hypothesis
that the 6-MP drug is inefficient, we would regard all patients as
belonging to a single group, instead of two. We could then implement
the latent diffusion model based on the stochastic perturbation of
the Weibull, but with just one diffusion process. Let M$_1$ denote the model
where all patients belong to a single group (corresponding to the
hypothesis $H_1$ of null efficacy of the 6-MP drug) and let M$_2$
denote the model considered above (corresponding to the hypothesis
$H_2$ of efficacy of the 6-MP drug). If the a priori probabilities of
hypotheses $H_1$ and $H_2$ are set equal to 0.5, the Bayes factor
\[
\mbox{BF}=\frac{\mbox{probability density of data under model
M$_1$}}{\mbox{probability density of data under model M$_2$}}
\]
gives the posterior odds in favor of $H_1$. As expected, the
computed Bayes factor ($\mathrm{BF}=9\times10^{-6}$) provides strong evidence
for the efficacy of the 6-MP drug.

\section{Latent diffusion models with covariates}
\label{sec:covariates}

Covariates can be included in the latent diffusion models described
in a very natural way, as directly influencing the underlying
diffusion. For instance, if $\mathbf{ Z}$ is a vector of $p$ covariates
measured at time 0, we can use the model based on the diffusion
satisfying the SDE
\begin{eqnarray}
\label{eq:SDE_covariates}
 \mathrm{d}X_t &=& \beta(X_t,\mathbf{ z},\theta)\,\mathrm{d}t+\sigma \,\mathrm{d}B_t,\qquad t\geq0,\nonumber
 \\[-8pt]\\[-8pt]
 X_0 &=& x_0(\mathbf{ z},\theta).\nonumber
\end{eqnarray}
In particular, following suggestions of
\citet{Aalen-Gjessing-2001} and\break \citet{Aalen-Borgan-Gjessing-2008} for
barrier hitting models, those covariates which represent measures of
how far the underlying process that leads to the event has
advanced (such as staging measures in cancer) may be taken to
influence the starting point of the diffusion. Those covariates
which instead represent causal influence on the development of the
process may be taken to influence the drift of the diffusion.

Let $\mathbf{ z}$ take values $\mathbf{ z}^{[1]},\ldots,\mathbf{ z}^{[q]}$.
Then, (\ref{eq:SDE_covariates}) gives $q$ different diffusions,
$X^{{[\mathbf{ z}=\mathbf{ z}^{[1]}}]}, \ldots,
X^{{[\mathbf{ z}=\mathbf{ z}^{[q]}}]}$, driven by the
same Brownian motion $B$, with
\begin{eqnarray*}
 \mathrm{d}X_t^{{[\mathbf{ z}=\mathbf{ z}^{[k]}}]}&=&\beta
\bigl(X_t^{{[\mathbf{ z}=\mathbf{ z}^{[k]}}]},
\mathbf{ z}^{[k]},\theta\bigr)\,\mathrm{d}t+\sigma \,\mathrm{d}B_t,\qquad t\geq0,
\\
 X_0 &=& x_0\bigl(\mathbf{ z}^{[k]}\bigr)
\end{eqnarray*}
for $k=1,\ldots,q$. Denote by $F_{X^{{[\mathbf{
z}=\mathbf{ z}^{[1]}}]},h}, \ldots, F_{X^{{[\mathbf{
z}=\mathbf{ z}^{[q]}}]},h}$ the relative random distributions, as in
(\ref{DDRM}). Moreover, denote by $Y_1^{{[\mathbf{
z}=\mathbf{ z}^{[k]}}]},\ldots,Y_{n_k}^{{[\mathbf{ z}=\mathbf{
z}^{[k]}}]}$ the survival times of the $n_k$ individuals having
covariates $\mathbf{ z}=\mathbf{ z}^{[k]}$ for $k=1,\ldots,q$. The
survival times $Y_1^{{[\mathbf{ z}=\mathbf{
z}^{[k]}}]},\ldots,Y_{n_k}^{{[\mathbf{ z}=\mathbf{
z}^{[k]}}]}$, conditionally on $F_{X^{{[\mathbf{
z}=\mathbf{ z}^{[k]}}]},h}$, are i.i.d. with common distribution
$F_{X^{{[\mathbf{ z}=\mathbf{ z}^{[k]}}]},h}$. Since the
$q$ diffusions are driven by the same Brownian motion, it is here
more natural to use the fully non-centered parametrization of the
model, described in Section \ref{pnc-section}. In particular, the
joint distribution of $Y_1^{{[\mathbf{ z}=\mathbf{
z}^{[1]}}]},\ldots,Y_{n_1}^{{[\mathbf{ z}=\mathbf{
z}^{[1]}}]}, \ldots, Y_1^{{[\mathbf{ z}=\mathbf{
z}^{[q]}}]},\ldots,Y_{n_q}^{{[\mathbf{ z}=\mathbf{
z}^{[q]}}]},$ given $B=b$ and $\Theta=\theta$, has density, with
respect to $\mathcal{L}^n$ (where $n=n_1+\cdots+n_q$), given by
\begin{eqnarray*}
&&l\bigl(y_1^{{[\mathbf{ z}=\mathbf{
z}^{[1]}}]},\ldots,y_{n_1}^{{[\mathbf{ z}=\mathbf{
z}^{[1]}}]}; \ldots;y_1^{{[\mathbf{ z}=\mathbf{
z}^{[q]}}]},\ldots,y_{n_q}^{{[\mathbf{ z}=\mathbf{
z}^{[q]}}]}\big\vert \theta,b_{[0,y_{[n]}]},\mathbf{ z}^{[1]},\ldots,\mathbf{
z}^{[q]}\bigr)
\\
&&\quad=
\prod_{k=1}^{q}l\bigl(y_1^{{[\mathbf{ z}=\mathbf{
z}^{[k]}}]},\ldots,y_{n_k}^{{[\mathbf{ z}=\mathbf{
z}^{[k]}}]}\big\vert \theta,b_{[0,y_{[n_k]}]},\mathbf{ z}^{[k]}\bigr),
\end{eqnarray*}
where $y_{[n]}:=\max\{y_1,\ldots,y_n\}$,
$y_{[n_k]}:=\max\{y_1^{{[\mathbf{ z}=\mathbf{
z}^{[k]}}]},\ldots,y_{n_k}^{{[\mathbf{ z}=\mathbf{
z}^{[k]}}]}\}$ and
\begin{eqnarray*}
l\bigl(y_1^{{[\mathbf{ z}=\mathbf{
z}^{[k]}}]},\ldots,y_{n_k}^{{[\mathbf{ z}=\mathbf{
z}^{[k]}}]}\big\vert \theta,b_{[0,y_{[n_k]}]},\mathbf{
z}^{[k]}\bigr)=l\bigl(y_1^{{[\mathbf{ z}=\mathbf{
z}^{[k]}}]},\ldots,y_{n_k}^{{[\mathbf{ z}=\mathbf{
z}^{[k]}}]}\big\vert x_{[0,y_{[n_k]}]}^{{[\mathbf{ z}=\mathbf{
z}^{[k]}}]}\bigr)
\end{eqnarray*}
is as in (\ref{likelyhood}). The joint posterior distribution of
$\Theta$ and $B$ has density, with respect to the product measure
$\mathcal{L}^d\otimes\mathbb{W}_{\sigma}$, given by
\begin{eqnarray}
\label{eq:post_covariates}
&&\pi\bigl(\theta,b\big\vert y_1^{{[\mathbf{ z}=\mathbf{
z}^{[1]}}]},\ldots,y_{n_1}^{{[\mathbf{ z}=\mathbf{
z}^{[1]}}]}; \ldots;y_1^{{[\mathbf{ z}=\mathbf{
z}^{[q]}}]},\ldots,y_{n_q}^{{[\mathbf{ z}=\mathbf{
z}^{[q]}}]}; \mathbf{ z}^{[1]},\ldots,\mathbf{ z}^{[q]}\bigr)\nonumber
\\[-8pt]\\[-8pt]
&&\quad =C
p_{\Theta}(\theta)\prod_{k=1}^{q}
l\bigl(y_1^{{[\mathbf{ z}=\mathbf{
z}^{[k]}}]},\ldots,y_{n_k}^{{[\mathbf{ z}=\mathbf{
z}^{[k]}}]}\big\vert \theta,b_{[0,y_{[n_k]}]},\mathbf{ z}^{[k]}\bigr).\nonumber
\end{eqnarray}

Note that this model is structurally different from the model for
multiple groups of observations described in Section
\ref{sec:multiple_groups} since the distributions of the survival
times are here linked at the level of the Brownian motion, allowing
a much stronger borrowing of strength for inference across
individuals who share a common value of even just one of the $p$
covariates.

As usual, we denote by $T$ the time horizon of interest, $T\geq
y_{[n]}$. The Hastings-within-Gibbs algorithm for sampling from
(\ref{eq:post_covariates}) alternates between
\begin{enumerate}[1.]
\item simulation of $\Theta$, conditional on the current path of
$B_{[0,y_{[n]}]}$, the
observations and the covariates;
\item simulation of $B_{[0,T]}$, conditional on the current value of
$\Theta$, the observations and the covariates.
\end{enumerate}
In particular, the update of the Brownian motion $B_{[0,T]}$ can be
carried out via the updating strategy based on overlapping Brownian
bridges, as described in Section \ref{pnc-section}.

\subsection{An illustrative application to a real-world data set with
covariates}\label{sec:appl_covariates}

In this section, we illustrate how to efficiently handle the model
with covariates via an application to a data set concerning 272
patients diagnosed with non-small cell lung cancer. The data set is
described in detail in \citet{Muers-Shevlin-Brown-1996}. Survival
times are measured in months from the time of diagnosis (with 17\%
of censoring) and some covariates are recorded at the time of
diagnosis. Just to give an illustration of the model, we shall
consider here two covariates: sex ($F=0$: male and $F=1$: female) and
hoarseness ($H=0$: absent and $H=1$: present). Using, for instance, the
model based on the stochastic perturbation around the Weibull, we
can include these covariates as follows:
\begin{eqnarray*}
\mathrm{d}X_t&=&\exp\{\theta_{10}+\theta_{11}F\} (\operatorname{sign}(X_t)) |X_t|^{\theta
_2}\, \mathrm{d}t + \sigma \,\mathrm{d}B_t, \qquad t\geq 0,
\\
X_0&=& \exp\{\theta_{00}+\theta_{01}F+\theta_{02}H\}.
\end{eqnarray*}
Note that, following the suggestion of
\citet{Aalen-Borgan-Gjessing-2008}, we have modeled the covariate
hoarseness, which only represents a
measure of how far the lung tumor has advanced, as influencing the
starting point of the diffusion; we have instead taken the covariate
sex to influence both the starting point and the drift of the
diffusion, in order to account for possible differences between
males and females, both in the hazards at time of diagnosis and in
the hazard dynamics. The covariate combinations determine four
different diffusions, $X^{{[F=0,H=0]}}$,
$X^{{[F=0,H=1]}}$,
$X^{{[F=1,H=0]}}$ and
$X^{{[F=1,H=1]}}$, driven by the same Brownian
motion. According to this model, the hazard at time 0 (the time of
diagnosis) of patients suffering from hoarseness is
$\exp\{\theta_{02}\}$ times that of patients not suffering from
hoarseness and the hazard at time 0 of female patients is
$\exp\{\theta_{01}\}$ times that of male patients; moreover,
$\exp\{\theta_{11}\}$ gives a measure of the different progression
rate of the cancer in female patients with respect to male
patients.

We express the data as fractions of a quadrennium and choose as
time horizon $T$ the maximum of the observations, corresponding to
about 37 months. In order to avoid dependencies among the
$(\theta_{00},\theta_{01},\theta_{02})$ parameters and among the
$(\theta_{10},\theta_{11})$ parameters, we reparametrize them in
terms of $(\eta_{00},\theta_{01},\theta_{02})$ and
$(\eta_{10},\theta_{11})$, with
$\theta_{00}=\eta_{00}-p_F\theta_{01}-p_H\theta_{02}$ and
$\theta_{10}=\eta_{10}-p_F\theta_{11}$, where we have denoted by
$p_F$ and $p_H$ the percentage of females patients and the
percentage of patients suffering from hoarseness, respectively. We
take all of the parameters to be a priori independent, with Gaussian
priors with mean $0$ and variance $5$ for all parameters except $\Theta
_2$, for which we use a uniform prior over $[0,1]$. Moreover, we set
$\sigma=8$. We then run the Hastings-within-Gibbs
algorithm based on the non-centered parametrization of the model.
The update of the parameters is performed via independence samplers
having proposal distributions equal to the priors. The algorithm is
run for 200\,000 iterations and the
first 2000 are discarded as burn-in.

\begin{figure}[t]

\includegraphics{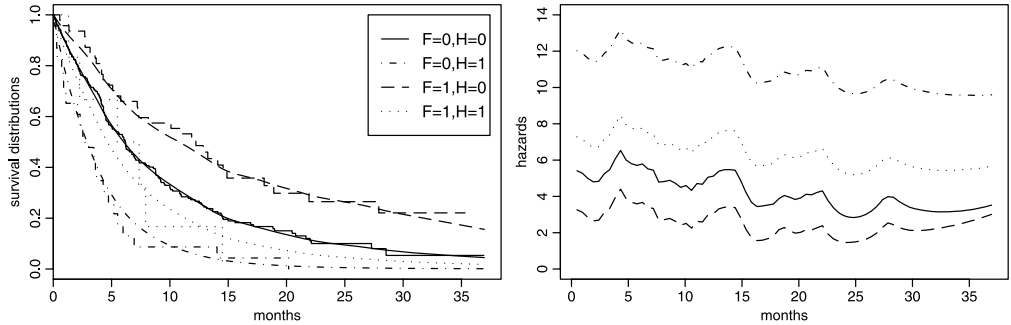}

\caption{Left:
posterior mean survival distributions, together with Kaplan--Meier
curves, for male patients without hoarseness at time of diagnosis
($F=0,H=0$, solid line), for male patients with hoarseness
($F=0,H=1$, dotted and dashed line), for female patients without
hoarseness ($F=1,H=0$, dashed line) and for female patients with
hoarseness ($F=1,H=1$, dotted line). Right: the same for posterior
mean hazard functions.} \label{fig:surv-haz-lungcancer-all}
\end{figure}

\begin{figure}[t]

\includegraphics{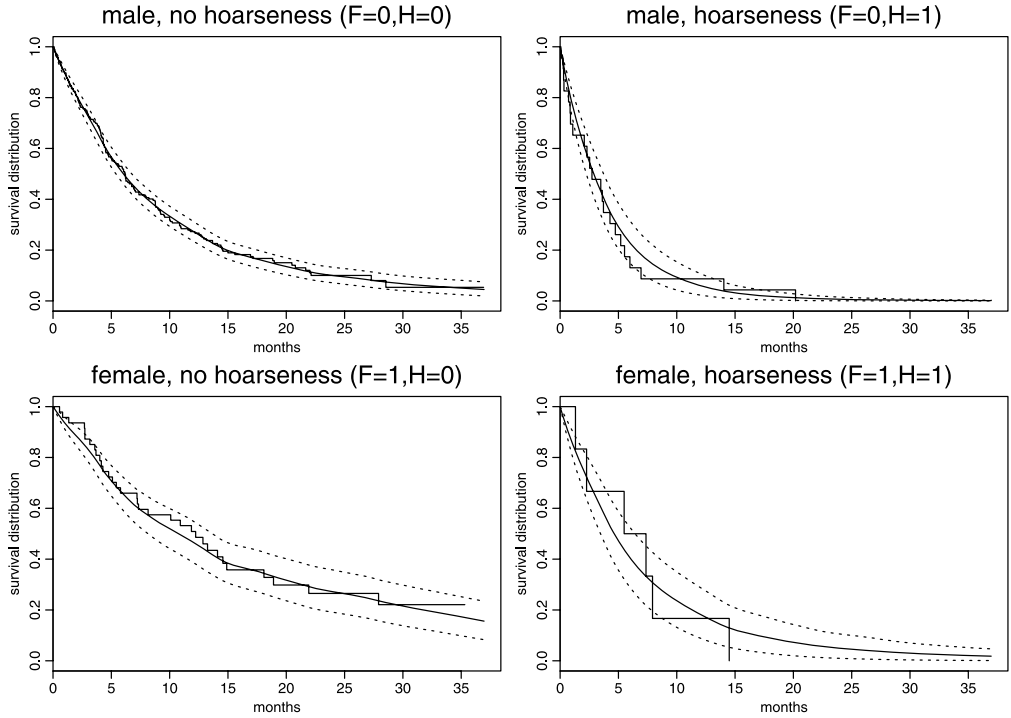}

\caption{Upper-left: posterior mean survival distribution and pointwise approximate
$90\%$ highest posterior bands, together with Kaplan--Meier curve,
for male patients without hoarseness. Upper-right: the same for male
patients with hoarseness. Lower-left: the same for female patients
without hoarseness. Lower-right: the same for female patients with
hoarseness.} \label{fig:surv-lungcancer-4plots}
\end{figure}

Figure \ref{fig:surv-haz-lungcancer-all} shows posterior mean
survival distributions, together with Kaplan--Meier curves, for male
patients without hoarseness at time of diagnosis ($F=0,H=0$, solid
line), for male patients with hoarseness ($F=0,H=1$, dotted and dashed
line), for female patients without hoarseness ($F=1,H=0$, dashed
line) and for female patients with hoarseness ($F=1,H=1$, dotted
line). The four survivals are also plotted separately in Figure
\ref{fig:surv-lungcancer-4plots} with $90\%$ highest posterior
bands. Figure
\ref{fig:surv-haz-lungcancer-all} also
displays the posterior mean hazard functions for the four covariate
combinations. In particular, the posterior mean hazard at time 0 of
patients suffering from hoarseness is 2.2 times bigger than that of
patients not suffering from hoarseness, whereas the hazard at time~0
of female patients is 0.6 times that of male patients.

Note that even though we have only
considered categorical covariates in this illustrative application,
quantitative covariates can
also be included in the model; however, it may be necessary to categorize
these covariates in order to have a sufficient number of
observations for each of the diffusion processes. This, of course,
requires larger data sets.


\section{Generalization to the case of unknown diffusion coefficient}
\label{sigma-unknown}

An important generalization of the models we have considered thus far
consists of considering diffusion processes with unknown
diffusion coefficient $\sigma$ since $\sigma$ describes a natural
measure of prior uncertainty. We briefly discuss how to deal with
this case.

Let $\Sigma$ be a real random variable. Given $\Theta=\theta$ and
$\Sigma=\sigma$, consider the scalar diffusion process $X$ solution
of the SDE (\ref{SDE}) and denote by $\mathbb{P}_{T,\theta,\sigma}$
the law of $X_{[0,T]}$. Let $p_{\Sigma}(\cdot)$ be the prior
density, with respect to $\mathcal{L}$, of $\Sigma$ (for
simplicity, we take $\Theta$ and $\Sigma$ to be stochastically
independent a priori). Let us consider, for instance, the centered
parametrization of the model. The joint posterior distribution of
$(\Theta,\Sigma,X_{[0,T]})$ has density, with respect to
$\mathcal{L}^{d+1}\otimes\mathbb{W}_{T,\sigma}$, given by
%
\begin{equation}
\label{posterior-sigma}
\pi\bigl(\theta,\sigma,x_{[0,T]}\vert y_1,\ldots,y_n\bigr)=C
p_{\Theta}(\theta) p_{\Sigma}(\sigma) g\bigl(x_{[0,T]}\vert \theta,
\sigma\bigr) l\bigl(y_1,\ldots,y_n\vert x_{[0,y_{[n]}]}\bigr),
\end{equation}
where $C$ is a normalizing constant and $g(x_{[0,T]}\vert\theta,
\sigma):=\frac{\mathrm{d}\mathbb{P}_{T,\theta,
\sigma}}{\mathrm{d}\mathbb{W}_{T,\sigma}}(x_{[0,T]})$ is given by Girsanov's
formula (\ref{girsanov}).

The
quadratic
variation of a diffusion processes, having diffusion coefficient
$\sigma$, satisfies
\begin{eqnarray*}
\lim_{m\to
\infty}\sum_{i=1}^{m}\bigl(X_{ti/m}-X_{t(i-1)/m}\bigr)^2=t\sigma^2,\qquad
\mathbb{W}_{T,\sigma}\mbox{-a.s.  for all }t.
\end{eqnarray*}
Therefore, the conditional distribution of $\Sigma$,
given the diffusion $X_{[0,T]}$, degenerates to a point mass and
$\Sigma$ is completely determined by the diffusion path. In
practice, we cannot simulate the diffusion path in continuous time,
but just at discrete time instants. In any case, the finer the discrete-time
approximation $\{X_{iT/m}\dvt i=1,\ldots,m\}$ of the diffusion
$X_{[0,T]}$, the stronger the dependence between
$\{X_{iT/m}\dvt i=1,\ldots,m\}$ and $\Sigma$. Consider the algorithm for
the simulation from (\ref{posterior-sigma}) that alternates
between:
\begin{enumerate}[1.]
\item simulation of
$\Theta$, conditional on the current value of $\Sigma$ and the
current path of
$X_{[0,T]}$;
\item simulation of $\Sigma$, conditional on the current value of
$\Theta$ and the current path of
$X_{[0,T]}$;
\item simulation of $X_{[0,T]}$, conditional on the observations
and the current values of $\Theta$ and $\Sigma$.
\end{enumerate}
The finer the approximation of the diffusion path, the worse the
convergence of the algorithm becomes. In the limiting case
$m=\infty$ (i.e., if the diffusion process could be simulated in
continuous time), this scheme would be reducible; see
\citet{Roberts-Stramer-2001}. An alternative way to see this problem
is to note that the collection of measures
$\{\mathbb{W}_{T,\sigma}\dvt \sigma\in\mathbb{R}\}$ are mutually
singular and, therefore, so are the measures
$\{\mathbb{P}_{T,\theta,\sigma}\dvt \sigma\in\mathbb{R}\}$.

In this case, the need for a different parametrization of the model
is thus compelling. Following \citet{Roberts-Stramer-2001}, we
parametrize the model in terms of $(\Theta,\Sigma,\dot{X})$, where
$\dot{X}_t=(X_t-X_0)/\Sigma$. By It\^{o}'s formula,
\begin{eqnarray*}
\mathrm{d}\dot{X}_t=\frac{\beta(\dot{X}_t,\Theta)}{\Sigma}\,\mathrm{d}t+\mathrm{d}B_t,\qquad
t\geq0,\qquad\dot{X}_0=0.
\end{eqnarray*}
The distribution of $\dot{X}_{[0,T]}$ depends on $\Sigma$, but any
realization of $\dot{X}_{[0,T]}$ contains only finite information
about $\Sigma$. Analogous reparametrizations are derived starting
from the ones described in Section \ref{pnc-section}.
MCMC algorithms based on these reparametrizations can be obtained as
simple modifications of the ones previously described.

\begin{figure}[t]

\includegraphics{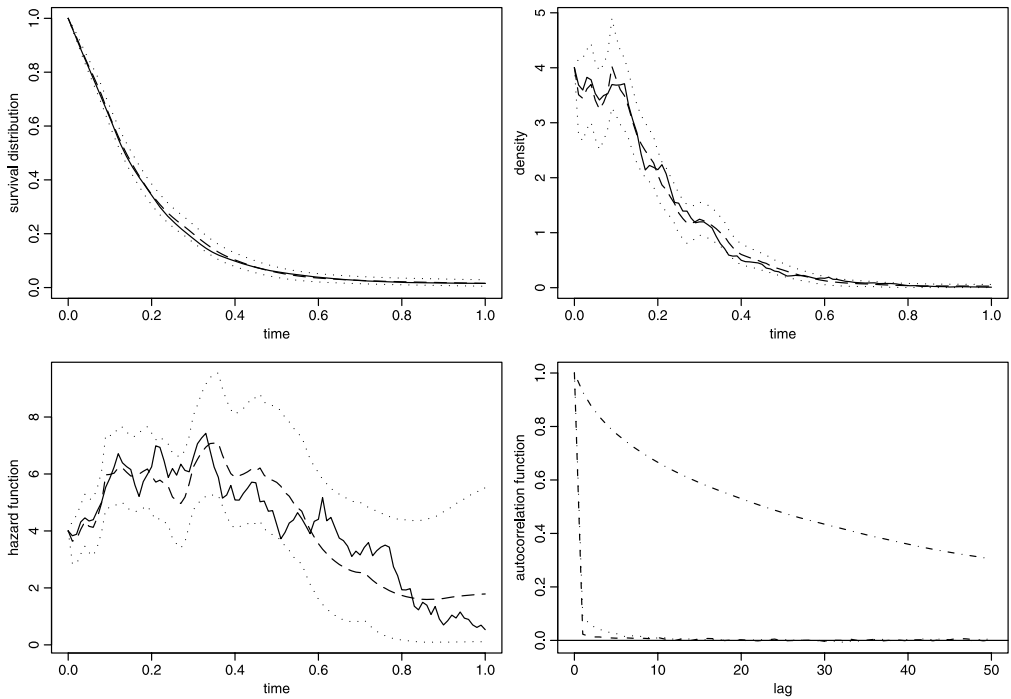}

 \caption{This corresponds to Figure
\protect\ref{fig:risultati}, but for the model with unknown diffusion
coefficient. The lower-right plot also displays the autocorrelation function
for $\sigma$ series (dotted-dashed line).}
\label{fig:risultati-c-sigma}
\end{figure}

Consider the
toy example described in Section \ref{simulation-p-e} and assume
the same model, but let the diffusion process have an unknown
diffusion coefficient. Let the prior for this coefficient be
exponential with mean 1. Figure \ref{fig:risultati-c-sigma} displays
the results obtained with the MCMC algorithm based on the
reparametrization $(\Theta,\Sigma,\dot{X})$. Specifications of the
algorithm are as in Section \ref{simulation-p-e}.
Note that the mixing for $\sigma$ is
slow relative to the very good mixing for $\theta_1$ and
$\theta_2$, but this does not prevent good estimates of the survival
distribution, density and hazard being obtained.
Slow
mixing for $\sigma$ could probably be improved by a further
reparametrization of the model.

Alternatively to the case of an unknown diffusion coefficient, it would
be possible to consider models based on diffusion processes having
$\sigma=1$, but with hazard function $h(\Gamma,X)$, where $\Gamma$
is a random parameter. A reparametrization of the
model would also be necessary in this case.

\section{Discussion}
\label{sec:discussion}

In this paper, we have described latent diffusion models for
survival analysis and have shown that these models can be
efficiently treated by means of MCMC techniques. We have dealt with
the case of multiple groups of observations, typical of clinical
trials, and we have shown how covariates can be efficiently included
in the models. We have outlined how, in the described framework, it is
possible to consider stochastic perturbations of common survival
models. In particular, we have used a stochastic perturbation of the
Weibull model in some illustrative applications to small data sets,
with multiple groups of observations and with covariates.
Applications to larger data sets, where the potential of a latent
diffusion model may be fully expressed, will be the object of future
work. All analyses presented are computationally feasible within R
(see \citet{manualR}).

Another generalization of the model we intend to explore regards
random probabilities based on jump diffusion processes. As observed
in Section \ref{prelim}, the cumulative hazard functions associated
with random probabilities based on diffusions are smooth, being the
integrals of continuous processes. By replacing the diffusion
process with a jump diffusion process, it would be possible to
capture sudden changes
in the behavior of cumulative hazards that might be due to some kind of
shock experienced by the
population. Hazards modeled through stochastic processes with jumps
have been studied, for instance, by \citet{Gjessing-Aalen-Hjort-2003}.

\section*{Acknowledgments}

We would like to thank Robin Henderson and Piercesare Secchi for
useful comments, and Omiros Papaspiliopoulos and Alexandros Beskos
for their help with programming. We are also grateful to the
Associate Editor and two anonymous referees for their constructive
comments. The second author acknowledges funding from the EC Marie Curie
Training Site Human Potential Program in order to visit the Department of
Mathematics and Statistics, Lancaster University and from the Centre
for Research in Statistical Methodology (CRiSM), University of
Warwick.

\printhistory

\end{document}